\documentclass[10pt,twosided]{article}
\usepackage[tbtags]{amsmath}
\usepackage{epsfig,amstext,amssymb,amsthm,latexsym}
\allowdisplaybreaks[4]
\usepackage{color}

\usepackage{amssymb,color}

\definecolor{c20}{rgb}{0.,0.7,0.}
\definecolor{c30}{rgb}{0.,0.,1.}
\definecolor{c40}{rgb}{1,0.1,0.7}
\definecolor{c50}{rgb}{1,0,0}
\definecolor{c60}{rgb}{1,0.9,0.1}
\def\Ke#1{\textcolor{c30}{#1}}
\def\Ke#1{#1}

\def\rH#1{\textcolor{c30}{#1}}
\def\rH#1{#1}
\def\cE#1{\textcolor{c30}{#1}}
\def\cF#1{\textcolor{c30}{#1}}
\def\ce#1{\textcolor{c20}{#1}}
\def\ce#1{#1}
\def\cE#1{#1}

\newcommand{\tb}[1]{{\textcolor{blue}{#1}}}
\def\tb#1{#1}

\def\cE#1{#1}
\def\cF#1{#1}
\def\cL#1{\textcolor{c30}{#1}}
\def\cL#1{#1}

\usepackage{amssymb,color}

\newcommand{\kb}[1]{\boldsymbol{#1}}
\newcommand{\vk}[1]{\kb{#1}}

\newcommand{\ve}{\varepsilon}

\newcommand{\abs}[1]{\lvert #1 \rvert}

\newcommand{\E}[1]{\mathbb{E}\left(#1\right)}
\newcommand{\VA}[1]{\mathbb{V}ar\left(#1\right)}
\newcommand{\CO}[1]{\mathbb{C}ov\left(#1\right)}

\newcommand{\pk}[1]{\mathbb{P} \left\{ #1 \right\} }

\newcommand{\pb}[1]{\mbox{\rm$\vk{P}$}\Bigl \{#1 \Bigr \}}

\newcommand{\todis}{\stackrel{d}{\to}}

\newcommand{\BQN}{\begin{eqnarray}}
\newcommand{\EQN}{\end{eqnarray}}
\newcommand{\BQNY}{\begin{eqnarray*}}
\newcommand{\EQNY}{\end{eqnarray*}}

\newcommand{\BS}{\begin{sat}}
\newcommand{\ES}{\end{sat}}
\newcommand{\BT}{\begin{theo}}
\newcommand{\ET}{\end{theo}}
\newcommand{\BK}{\begin{korr}}
\newcommand{\EK}{\end{korr}}

\newcommand{\BD}{\begin{de}}
\newcommand{\ED}{\end{de}}

\newcommand{\BIT}{\begin{itemize}}
\newcommand{\EIT}{\end{itemize}}
\newcommand{\BDI}{\begin{description}}
\newcommand{\EDI}{\end{description}}

\newcommand{\BRM}{\begin{remarks}}
\newcommand{\ERM}{\end{remarks}}

\newcommand{\BEL}{\begin{lem}}
\newcommand{\EEL}{\end{lem}}

\newtheorem{theo}{Theorem}[section]
\newtheorem{sat}[theo]{Proposition}
\newtheorem{de}[theo]{Definition}
\newtheorem{lem}[theo]{Lemma}

\newtheorem{example}[theo]{Example}
\newtheorem{korr}[theo]{Corollary}

\newtheorem{remarks}[theo]{Remarks}

\newcommand{\netheo}[1]{{Theorem \ref{#1}}}
\newcommand{\nekorr}[1]{{Corollary \ref{#1}}}

\newcommand{\prooftheo}[1]{ \textsc{Proof of Theorem} \ref{#1} }

\newcommand{\proofkorr}[1]{\textsc{Proof of Corollary} \ref{#1}}

\newcommand{\COM}[1]{}

\newcommand{\QED}{\hfill $\Box$}

\date{}

\def\IF{\infty}
\def\rw{\rightarrow}

\def\DT{\widetilde{\delta}(t)}
\def\DTT{\widetilde{\delta}(T)}
\def\oo{(1+o(1))}
\def\vn{\varepsilon}
\def\H{\mathcal{P}}

\begin{document}
\title{\bf \Large Gaussian Risk Models \cF{with Financial Constraints}} 
\author{
Krzysztof D\c{e}bicki\footnote{Mathematical Institute, University of Wroc\l aw, pl. Grunwaldzki 2/4, 50-384 Wroc\l aw, Poland},
 \ Enkelejd Hashorva$^\dag$, 
Lanpeng Ji\footnote{Department of Actuarial Science, University of Lausanne, UNIL-Dorigny 1015 Lausanne, Switzerland}
}\bigskip

 \maketitle

\begin{quote}
{\bf Abstract:} In this paper we investigate Gaussian risk models which include financial elements such as
inflation and interest rates. For some general models for inflation and interest rates, we \Ke{obtain}
\Ke{an} asymptotic expansion of the finite-time ruin probability for Gaussian risk models. Furthermore,
\rH{we derive an approximation of the conditional ruin time by 
 an} exponential random variable as the initial capital tends to infinity.

{\bf Key Words:} Finite-time ruin probability; \rH{conditional ruin time}; exponential approximation; Gaussian risk process; inflation; interest.

{\bf AMS Classification:} Primary 91B30; Secondary 60G15, 60G70.

\end{quote}

\section{Introduction }
A central topic in the actuarial literature, inspired by the early
contributions of Lundberg (1903) and Cram\'{e}r (1930), is
the computation of the ruin probability over both finite-time and
infinite-time horizon; see e.g., Rolski et al. (1999), Mikosch
(2008), Asmussen and Albrecher (2010) and \rH{the} references therein. As
mentioned in Mikosch (2008) calculation of \rH{the} ruin probability is
considered as the \tb{"jewel"} of the actuarial mathematics.\\
\rH{In fact},
exact formulas for both finite-time and infinite-time ruin
probability are known only for few special models. Therefore,
asymptotic methods have been developed to derive expansions of the
ruin probability as the initial capital/reserve increases to
infinity. Following Chapter 11.4 in Rolski et al.\ (1999) the risk
reserve process of an insurance company can be modelled by a
stochastic process $\{\tilde U(t),t\ge 0\}$ given as
\BQN\label{U0} \tilde U(t)=u+ct-\int_0^tZ(s)\, ds,\ \ \ t\ge0,
\EQN where $u\ge0$ is the initial reserve, $c>0$ is the rate of
premium received by the insurance company,
and $\{Z(t),t\ge0\}$ is a centered Gaussian process with almost surely continuous sample paths; 
the process $\{Z(t),t\ge0\}$ is \rH{frequently referred to as} the loss rate of the insurance company.
Under the assumption that  $\{Z(t),t\ge0\}$ is stationary the asymptotics of the infinite-time ruin probability of the process \eqref{U0}
 defined by
$$
\psi_{\IF}(u)=\pk{\inf_{t\in[0,\IF)} \tilde U(t)<0},\ \ \ u\ge0
$$
has been investigated in H\"usler and Piterbarg (2004), D\c{e}bicki (2002), Dieker (2005) and  Kobelkov (2005); see also H\"usler and Piterbarg (1999) \rH{and Hashorva et al.\ (2013).}
Therein the exact speed of convergence to 0 of $\psi_{\IF}(u)$ as $u\rw\IF$ was \rH{dealt with}.\\
In order to account for \rH{the} financial nature of the risks and thus for the time-value of the money \rH{as well as} other important
economic factors,
in this paper we shall consider a more general risk process which includes inflation/deflation effects and
interest rates (cf. Chapter 11.4 in Rolski et al. (1999)). \Ke{Essentially,} in case of inflation, a monetary unit at time 0 has the value $e^{-\delta_1(t)}$ at time $t$, where $\delta_1(t), t\ge0$ is a positive function with $\delta_1(0)=0$. In case of interest, a monetary unit invested at time 0 has the value $e^{\delta_2(t)}$ at time $t$, where $\delta_2(t), t\ge0$ is another positive function with $\delta_2(0)=0$.\\
Assuming first that the premium rate and the loss rate have to be adjusted for inflation, we
arrive at the following risk reserve process
\BQN
u+c\int_0^te^{\delta_1(s)}ds - \int_0^te^{\delta_1(s)}Z(s)ds,\ \ \ t\ge0.\nonumber
\EQN
Since the insurance company invests the surplus and thus accounting for investment effects  the resulting risk reserve
process is
\BQN\label{U1}
U(t)=e^{\delta_2(t)}\left(u+c\int_0^te^{\delta_1(s)-\delta_2(s)}
ds-\int_0^te^{\delta_1(s)-\delta_2(s)}
Z(s)ds\right),\ \ \ t\ge0.
\EQN
We shall refer to $\{U(t),t\ge 0\}$ as the {\it risk reserve process in an economic environment};  see Chapter 11.4 in Rolski et al.\ (1999) for a detailed discussion on the effects of financial factors on the risk reserve processes.

In the case that $\delta_1(t)=0,  \delta_2(t)=\delta t, t\ge0,$ with $\delta>0$, the random process $\{U(t), t\ge 0\}$ reduces to a {\it risk reserve process with constant force of interest}.
For a class of stationary Gaussian processes
$\{Z(t),t\ge0\}$ with twice differentiable covariance function,
the exact asymptotics of the infinite-time ruin probability for
the {risk reserve process with constant force of interest} was
obtained in He and Hu (2007). Since therein the authors considered
 only smooth Gaussian process, the method of proof relied
on the well-known Rice method; see e.g., Piterbarg (1996).

\cL{Let $T$ be any positive constant.}
The principal goal of this contribution is the derivation of the exact asymptotics
of the finite-time ruin probability of the risk reserve process $U$ given by
\BQN\label{Intruin1}
\psi_T(u)&:=&\pk{\inf_{t\in[0,T]}U(t)<0}\nonumber\\
&=&\pk{\sup_{t\in[0,T]}\biggl(\int_0^t e^{- \delta(s)} Z(s)\, ds-c\int_0^t e^{-\delta(s)}ds\biggr)>u}
\EQN
as $u\rw\IF$, where $\{Z(t),t\ge0\}$ is a general centered Gaussian process with almost surely continuous sample paths  and
$\delta(t)=\delta_2(t)-\delta_1(t), t\ge0$ is some measurable real-valued function satisfying $\delta(0)=0$. Note in passing that $\delta(t)>0$ means that the interest contributes
more to the risk reserve process than the inflation at time $t$, and vice versa.\\
\tb{
In \netheo{ThmM} \rH{below we shall} show that $\psi_T(u)$
has asymptotically, as $u\to\infty$, (non-standard) normal distribution. This emphases the
qualitative difference between asymptotics in finite- and infinite-time horizon scenario;
see He and Hu (2007).
}

\tb{A related, interesting} and vastly analyzed quantity is the time of ruin which in our model is \Ke{defined as}
\begin{eqnarray}
\tau(u)=\inf\{t\ge0: U(t)<0\}, \quad \Ke{u\ge 0}.\label{tau.def}
\end{eqnarray}

Using that
$\pk{\tau(u)<T}=\pk{\inf_{t\in[0,T]}U(t)<0}$,
investigation  of distributional properties of the time of ruin under the condition that
ruin occurs in \cL{a} certain time period \rH{has}  attracted substantial attention; see e.g.,
the seminal contribution Segerdahl (1955) and the monographs Embrechts et al. (1997) and Asmussen and Albrecher (2010).
 Recent results for \tb{infinite-time} Gaussian and L\'{e}vy risk models are derived in H\"{u}sler (2006), H\"{u}sler and Piterbarg (2008), H\"{u}sler and Zhang (2008),
  Griffin and Maller (2012), Griffin (2013) \rH{ and Hashorva and Ji (2013)}.

\newcommand{\equaldis}{\stackrel{d}{=}}

\tb{
In Theorem \ref{Thm2} we derive a novel result, which shows that as $u\to \IF$,
the sequence of random variables $\{\xi_u,u>0\}$, defined (on the same probability space) by
\BQN\label{xiu}
 \xi_u \equaldis  u^2(T-\tau(u))\Bigl \lvert (\tau(u)<T)
\EQN
converges in distribution
to an exponential random variable
(here  $\equaldis$ stands for the equality of the distribution functions).
This, somewhat surprising result, contrasts with
the infinite-time case analyzed by
H\"{u}sler and Piterbarg (2008) \rH{and Hashorva and Ji (2013)}, where the \rH{limiting random variable} is normally distributed.
}

\COM{we obtain a conditional limit theorem for $\tau(u)$ as $u\rw\IF$,
given that the ruin occurs before time $T$, i.e.,
\begin{eqnarray}
\pk{u^2(T-\tau(u))\le x|\tau(u)<T}
\label{lim0}
\end{eqnarray}
as $u\to\infty$.
For the infinite-time horizon  \ce{derived} the approximation of the conditional ruin time
by a normal random variable; our result is different and somewhat surprising showing that the correct approximation is achieved by an exponential distribution.
}
Organization of the paper: The main results concerning the finite-time ruin probability and the approximation of $\xi_u$ are displayed in Section 2, whereas the proofs are relegated to Section 3. \rH{We conclude this contribution with a short Appendix.}

\def\GuT{\frac{u+c\widetilde{\delta}(T)}{\sigma(T)}}
\section{Main Results}
Let the loss rate of the insurance company $\{Z(t),t\ge 0\}$ be modelled by a centered Gaussian
process with \rH{almost surely} continuous sample paths and covariance function $\CO{Z(s),Z(t)}=R(s,t)$.
As mentioned in the Introduction we shall require that $\delta(0)=0$. For notational simplicity we shall define below
\BQN\label{notation}
&&Y(t):= \int_0^t e^{- \delta(s)} Z(s) ds,\ \  \sigma^2(t):=\mathbb{V}ar(Y(t)),\ \  \DT:=\int_0^t e^{-\delta(s)}ds,\quad t\in [0,T]. 
\EQN
\tb{In what follows let $\sigma'(t)$
be the derivative of $\sigma(t)$, \rH{and let}
$\Psi$ denote the survival function of a $N(0,1)$ random variable.} We are interested in the asymptotic behavior of (\ref{Intruin1}) as the initial reserve $u$ tends to infinity,  i.e., we shall investigate
the asymptotics of
\[
\psi_T(u)=\pk{\sup_{t\in[0,T]}\biggl( Y(t)-c\DT\biggr)>u}
\]
as $u\to\infty$. In our first result below we derive an asymptotic expansion of $\psi_T(u)$ in terms of $c,\sigma(T),\widetilde{\delta}(T)$.

\BT\label{ThmM}
Let $\{Z(t), t\ge0\}$ be a centered Gaussian process with \rH{almost surely} continuous sample paths and covariance function
$R(s,t), s,t\ge0$. Further let  $\delta(t), t\ge0$, be some measurable function with $\delta(0)=0$.
If $\sigma(t)$ attains its maximum  over $[0,T]$ at \Ke{the} unique point $t=T$  and $\sigma'(T)>0$, 
then
\BQN\label{main}
\psi_T(u)\ =  \tb{\pk{\mathcal{N}>(u+ c \DTT)/\sigma(T)} \oo}
\label{IntruinTh1}
\EQN
holds as $u\rightarrow \infty$, with $\mathcal{N} $ a $N(0,1)$ random variable. 
\ET

\cF{The following result is an immediate consequence of \netheo{ThmM}.}
\begin{korr}\label{korr1}
Let $\{Z(t), t\ge0\}$ and $\delta(t), t\ge0$ be given as in \netheo{ThmM}.
If $R(s,t)>0$ for any $s,t\in[0,T]$, then \eqref{IntruinTh1} is satisfied.
\COM{\BQNY
\psi_T(u)\ =\frac{\sigma(T)}{\sqrt{2\pi}}u^{-1}\exp\left(-\frac{(u+c\cL{\DTT})^2}{2\sigma^2(T)}\right)\oo,
\EQNY
as $u\rightarrow \infty$.}
\end{korr}

\begin{remarks}
a) \tb{It follows from the proof of \netheo{ThmM} that  \eqref{main} still holds if
$Y(t):= \int_0^t e^{- \delta_1(s)} Z(s) ds$ and
$\DT:=\int_0^t e^{-\delta_2(s)}ds$
in (\ref{notation}).}
\\
b) In the asymptotic behavior of $\psi_T(u)$ the positive constant $\sigma'(T)$ does not appear. It appears however explicitly in the the approximation of the conditional ruin time as shown in our second theorem below.
\end{remarks}

Along with the analysis of the ruin probability  in risk theory an important theoretical topic is the behavior of the ruin time.
For infinite-time horizon results in this direction are well-known;  see e.g., Asmussen and Albrecher  (2010),
H\"{u}sler and Piterbarg (2008) and \rH{Hashorva and Ji (2013)} for the normal approximation of the conditional distribution of \Ke{the} ruin time $\tau(u)$ given that $\tau(u)<\IF$.\\
In our second result below we show that (appropriately rescaled) ruin time
$\tau(u)$ conditioned that $\tau(u)<T$ is asymptotically, as $u\to\infty$, exponentially
distributed with parameter
$\sigma'(T)/(\sigma(T))^3$.

\BT\label{Thm2}
Under the conditions of \netheo{ThmM}, we have
\BQN\label{tau}
\lim_{u\to\infty}
\pb{u^2(T-\tau(u))\le x \Bigl|\tau(u)<T}
=
1 - \exp\left(     -  \frac{\sigma'(T)}{\sigma^3(T)}x    \right),\ \ x\ge0.
\EQN
\ET
Note in passing that \eqref{tau} means the convergence in distribution
\BQN \label{meanE}
\xi_u \todis \xi, \quad u\to \IF,
\EQN
where $\xi$ is exponentially distributed \Ke{such that}
$$
 e_T:=\E{\xi}=  \frac{\sigma^3(T)}{\sigma'(T)}>0.
$$

We present next three illustrating examples.

\begin{example} Let $\{Z(t), t\ge0\}$ be an Ornstein-Uhlenbeck process with parameter $\lambda>0$, i.e., \Ke{$Z$ is} a stationary process with covariance function
$R(s,t)=\exp(-\lambda\abs{s-t})$. If $\delta(t)=\delta t, t\ge0$ with $\delta \in (0, \lambda)$, then
\BQNY
\DT=\frac{1}{\delta}\left(1-e^{-\delta t}\right),\ \ \ \ \sigma^2(t)=\frac{1}{(\lambda-\delta)\delta}\left(1-e^{-2\delta t}\right)-\frac{2}{\lambda^2-\delta^2}\left(1-e^{-(\lambda+\delta)t}\right), \ t\in[0,T].
\EQNY
Therefore, \tb{from \nekorr{korr1}, we obtain that}
\BQNY
\psi_T(u)\
&=&\frac{1}{u\sqrt{2 \pi}}
\sqrt{\frac{1}{(\lambda-\delta)\delta}\left(1-e^{-2\delta T}\right)-\frac{2}{\lambda^2-\delta^2}\left(1-e^{-(\lambda+\delta)T}\right)}
\\
&&\times
\exp\left(-
\frac{(u+ c/\delta \left(1-e^{-\delta T}\right))^2}{2\left(\frac{1}{(\lambda-\delta)\delta}\left(1-e^{-2\delta T}\right)-\frac{2}{\lambda^2-\delta^2}\left(1-e^{-(\lambda+\delta)T}\right)\right)}\right)
\oo
\EQNY
as $u\rightarrow \infty$. Furthermore, in view of \netheo{Thm2} the convergence in \eqref{meanE} holds with

\BQNY
e_T= \frac{((\lambda+\delta)(1-e^{-2\delta T}) -2\delta(1-e^{-(\lambda+\delta)T}))^2}    {(\lambda-\delta)(\lambda+\delta)^2\delta^2e^{-\delta T}(e^{-\delta T}-e^{-\lambda T}) }.
\EQNY
\end{example}

\begin{example} Let $\{Z(t), t\ge0\}$ be a Slepian process, i.e.,
$$
Z(t)=B(t+1)-B(t),\ \ \ t\ge0,
$$
with $B$ a standard Brownian motion. For this model we have
$R(s,t)=\max(1-\abs{s-t},0)$. If further  $\delta(t)=\delta t, t\ge0$ with $\delta\neq0$, then
\BQNY
\DT=\frac{1}{\delta}\left(1-e^{-\delta t}\right),\ \ \ \ \sigma^2(t)=\frac{1}{\delta^2}-\frac{1}{\delta^3}-\frac{2}{\delta^2}e^{-\delta t}+\frac{2t}{\delta^2}e^{-\delta t}
+\frac{\delta+1}{\delta^3}e^{-2\delta t}, \ t\in[0,1].
\EQNY
Consequently, \nekorr{korr1} implies, as $u\to \IF$
\BQNY
\psi_1(u)\ &=&
\sqrt{\frac{
\delta-1
+(\delta+1) e^{-2\delta }}{2\pi\delta^3}}u^{-1}
\exp\left(-
\frac{\delta\left(\delta u+ c  \left(1-e^{-\delta }\right)\right)^2}
     {2(\delta-1
+(\delta+1) e^{-2\delta })}
\right)\oo.
\EQNY
\Ke{Further} by \netheo{Thm2} the convergence in \eqref{meanE} holds with
\BQNY
e_1= \frac{(\delta-1+ (\delta+1) e^{-2\delta} )^2}{\delta^4e^{-\delta} -\delta^4(\delta+1)e^{-2\delta} }.
\EQNY
\end{example}

\begin{example}
Let $\{Z(t), t\ge0\}$ be a standard Brownian motion
and assume that $\delta(t)= t^2/2, t\ge0$. Since
$R(s,t)=\min(s,t)$ we obtain
\BQNY
\DT=\sqrt{2\pi}(1/2-\Psi(t)),\ \ \ \ \sigma^2(t)=(\sqrt{2}-1)\sqrt{\pi}-2\sqrt{2\pi}\Psi(t)+2\sqrt{\pi}\Psi(\sqrt{2}t), \ t\in[0,T].
\EQNY
Applying once again \nekorr{korr1} we obtain
\BQNY
\psi_T(u)\ &=&
\sqrt{\frac{(\sqrt{2}-1)\sqrt{\pi}-2\sqrt{2\pi}\Psi(T)+2\sqrt{\pi}\Psi(\sqrt{2}T)}{2\pi}}u^{-1}\\
&&\times
\exp\left(
-\frac{\left(u+\sqrt{2\pi}c(1/2-\Psi(T))\right)^2}{2\left((\sqrt{2}-1)\sqrt{\pi}-2\sqrt{2\pi}\Psi(T)+2\sqrt{\pi}\Psi(\sqrt{2}T)\right) }
\right)\oo
\EQNY
as $u\rightarrow \infty$. Finally, by \netheo{Thm2} the convergence in distribution in \eqref{meanE} holds with
\BQNY
e_T= \frac{\left((\sqrt{2}-1)\sqrt{\pi}-2\sqrt{2\pi}\Psi(T)+2\sqrt{\pi}\Psi(\sqrt{2}T)\right)^2}{\sqrt{2\pi}(\varphi(T)-\varphi(\sqrt{2}T))  },
\EQNY
where $\varphi=-\Psi'$ is the \tb{density function of \rH{$N(0,1)$} random variable}.
\end{example}

\section{Proofs}
Before presenting proofs of Theorems \ref{ThmM} and \ref{Thm2}, we introduce some notation.
Let
$g_u(t)= \frac{u+c\tilde\delta (t)}{\sigma(t)}$ and  define 
\BQNY
X_u(t)&:=&\frac{Y(t)}{\sigma(t)}\frac{g_u(T)}{g_u(t)},
\quad \sigma_{X_u}^2(t):=\VA{X_u(t)},\\
r_{X_u}(s,t)&:=&\CO{\frac{X_u(s)}{\sigma_{X_u}(s)},\frac{X_u(t)}{\sigma_{X_u}(t)}}
=\CO{\frac{Y(s)}{\sigma(s)},\frac{Y(t)}{\sigma(t)}}.
\EQNY
Then, we can reformulate
 (\ref{Intruin1}) for all large $u$ as
\BQN
\psi_T(u)&=&
\pk{\underset{t\in[0,T]}\sup\left(\frac{Y(t)}{\sigma(t)}\frac{g_u(T)}{g_u(t)}\right)\ >\ g_u(T)}\notag\\
&=&
\pk{\underset{t\in[0,T]}\sup X_u(t)\ >\ g_u(T)}
,\label{intruing} \ \ u\ge0.
\EQN

\prooftheo{ThmM}
\cE{We shall derive first a lower bound for $\psi_T(u)$.} It follows from (\ref{intruing}) that
\BQN
\psi_T(u)
&\ge&\pk{\frac{Y(T)}{\sigma(T)}\ >\ g_u(T)}\nonumber\\
& =& \Psi(g_u(T))\nonumber\\
&=& \frac{\sigma(T)}{\sqrt{2\pi}}u^{-1}\exp\left(-\frac{(u+c\DTT)^2}{2\sigma^2(T)}\right)\oo
\label{intruinlow}
\EQN
as $u\rw\IF$. Next, we derive the upper bound.
\cE{Since}  $R(s,t)=R(t,s)$ for any $s,t\in[0,T]$,
we have
\BQN
\label{intruinvar}
\sigma^2(t):=\VA{Y(t)} 
=2\int_0^t\int_0^we^{-\delta(v)-\delta(w)}R(v,w)\, dvdw.
\EQN
Further, since by the assumption the function $\sigma(t)$ attains its unique maximum
over $[0,T]$ at $t=T$ and that $\sigma'(T)>0$,
there exists some $\theta_1\in(0,T)$ such that \cL{$\sigma(t)$ is strictly increasing on $[\theta_1,T]$ and}
\BQN \label{intruinsig}
\inf_{ t\in[\theta_1,T]}\sigma'(t)>0 
\EQN
implying that for $u$ sufficiently large
\BQNY
\sigma'_{X_u}(t)=\frac{\sigma'(t)}{\sigma(T)}\frac{u+c\DTT}{u+c\DT}-\frac{ce^{\delta(t)}\sigma(t)(u+c\DTT)}{(u+c\DT)^2\sigma(T)}>0
\EQNY
for all $ t\in[\theta_1,T]$.
Hence,  for sufficiently large $u$, $\sigma_{X_u}(t)$ is strictly increasing on $[\theta_1,T]$.
Furthermore, since
\BQNY
1-\sigma_{X_u}(t)&=&1-\frac{g_u(T)}{g_u(t)}\nonumber\\
&=&\frac{(\sigma(T)-\sigma(t))\left(u+c\DT \right)-c\sigma(t)(\DTT-\DT)}{\sigma(T)\left(u+c\DT\right)},
\EQNY
then by the definitions of $\DT$ and $\sigma(t)$ for any $\vn_1>0$
there exist some constants $K>0$ and $\theta_2\in(0,T)$ such that
\BQNY
\DTT-\DT&\le& K\ (T-t),\nonumber\\
(1-\vn_1)\sigma'(T)(T-t)&\le&\sigma(T)-\sigma(t)\le (1+\vn_1)\sigma'(T)(T-t)
\EQNY
\cL{are} valid for all $t\in[\theta_2,T]$. Therefore, we conclude that for  $u$  sufficiently large
\BQN\label{DDt}
(1-\vn_1)^2 \frac{\sigma'(T)}{\sigma(T)}(T-t) \le1-\sigma_{X_u}(t)\le (1+\vn_1)\frac{\sigma'(T)}{\sigma(T)}(T-t)
\EQN
for $t\in[\theta_2,T]$. For any $s<t$ we have 
\BQNY
1-r_{X_u}(s,t)&=&1-\mathbb{C}ov\left(\frac{Y(s)}{\sigma(s)},\frac{Y(t)}{\sigma(t)}\right)\nonumber\\
&=& \frac{\VA{Y(t)-Y(s)}-(\sigma(t)-\sigma(s))^2}{2\sigma(s)\sigma(t)}\nonumber\\
&\le& \frac{\VA{Y(t)-Y(s)}}{2\sigma(s)\sigma(t)}\nonumber\\
&=&\frac{\int_s^t\int_s^t  \cL{R(v,w)}e^{-\delta(w)-\delta(v)} dwdv}
     {2\sigma(s)\sigma(t)}.
\EQNY
The above implies that for sufficiently large $u$ and $s,t\in[\theta_2,T]$
\BQN\label{intruincovX}
1-r_{X_u}(s,t)&\le&C(t-s)^2,
\EQN
where
$C=\max_{w,v\in [\theta_2,T]}\frac{\cL{\abs{R(v,w)}}e^{-\delta(w)-\delta(v)}}{2 \sigma^2(w)}$ .
Consequently, in the light  of \eqref{DDt} and \eqref{intruincovX}, for any $\vn>0$ sufficiently small,
we have for some $\theta_0\in(\max(\theta_1,\theta_2),T)$ 
\BQNY
\sigma_{X_u}(t)\le\frac{1}{1+(1-\ve)(1-\vn_1)^2\frac{\sigma'(T)}{\sigma(T)}(T-t)}
\EQNY
and
\BQNY
\ r_{X_u}(s,t)\ge e^{-(1+\ve)C (t-s)^2}
\EQNY
for all $s,t\in[\theta_0,T]$. Next, define a centered Gaussian process $\{Y_\ve(t),t\ge0\}$ as
$$
Y_\ve(t)=\frac{\xi_\ve(t)}{1+(1-\ve)(1-\vn_1)^2\frac{\sigma'(T)}{\sigma(T)}(T-t)},
$$
where $\{\xi_\ve(t),t\ge0\}$ is a centered stationary Gaussian process with covariance function
$\CO{\xi_\ve(t),\xi_\ve(s)}=e^{-(1+\ve)C (t-s)^2}$. 
In view of Slepian Lemma (cf. Adler and Taylor (2007) or Berman (1992)) we obtain
\BQN
\pk{\underset{t\in[\theta_0,T]}\sup\left(Y(t)-c\DT\right)\ >\ u }& =& \pk{\underset{t\in[\theta_0,T]}\sup X_u(t)\ >\ g_u(T)}\nonumber\\
&\le&\pk{\underset{t\in[\theta_0,T]}\sup\left( \frac{X_u(t)/\sigma_{X_u}(t)}{1+(1-\ve)(1-\vn_1)^2\frac{\sigma'(T)}{\sigma(T)}(T-t)}\right)\ >\ g_u(T)}\nonumber\\
&\le&\pk{\underset{t\in[\theta_0,T]}\sup Y_\ve(t)\ >\ g_u(T)}\nonumber \\
&=& \Psi(g_u(T))\oo\label{intruinslep}
\EQN
as $u\rw\IF$, where the last asymptotic equivalence follows from iii) of \netheo{ThmPiter} in Appendix. Moreover since for $u$ sufficiently large there exists some $\lambda\in(0,1)$ such that
\BQNY
\underset{t\in[0,\theta_0]}\sup\sigma_{X_u}(t)\ \le\ \underset{t\in[0,\theta_0]}\sup\frac{(1+\lambda)\sigma(t)}{\sigma(T)}\ \le\ \frac{(1+\lambda)\sigma(\theta_0)}{\sigma(T)}\ <\ 1
\EQNY
 and
\BQNY
\pk{\underset{t\in[0,\theta_0]}\sup X_u(t)\ >a}\le\pk{\underset{t\in[0,\theta_0]}\sup\frac{2Y(t)}{\sigma(T)}>a}\le\frac{1}{2}
\EQNY
for some positive number $a$,
 we get from Borell inequality (e.g., Piterbarg (1996)) that, for $u$ sufficiently large
\BQN
\pk{\underset{t\in[0,\theta_0]}\sup\left(Y(t)-c\DT\right)\ >\ u }& =& \pk{\underset{t\in[0,\theta_0]}\sup X_u(t)\ >\ g_u(T)}\nonumber\\
&\le& 2\Psi\left(\frac{(g_u(T)-a)\sigma(T)}{(1+\lambda)\sigma(\theta_0)}\right)
=o(\Psi(g_u(T)))
\label{intruinborell}
\EQN
as $u\to\infty$. Combining \eqref{intruinslep} and \eqref{intruinborell}, we conclude that
\BQNY
\psi_T(u)& \le&  \pk{\underset{t\in[0,\theta_0]}\sup X_u(t)\ >\ g_u(T)}+\pk{\underset{t\in[\theta_0,T]}\sup X_u(t)\ >\ g_u(T)}\nonumber\\
&=& \Psi(g_u(T))\oo\nonumber\\
&=& \frac{\sigma(T)}{\sqrt{2\pi}}u^{-1}\exp\left(-\frac{(u+c\DTT)^2}{2\sigma^2(T)}\right)\oo
\EQNY
as $u\rightarrow\infty$, \rH{which together with \eqref{intruinlow} establishes the proof.} \QED

\proofkorr{korr1} Since  $R(s,t)>0$ for any $s,t\in[0,T]$ \rH{it follows} from \eqref{intruinvar} that $\sigma(t)$ attains its unique maximum  over $[0,T]$ at $t=T$  and $\sigma'(T)>0$. Therefore, the claim follows immediately from \netheo{ThmM}. \QED

\prooftheo{Thm2} In the following we \rH{shall} use the same notation as in the proof of \netheo{ThmM}. First note that for any $x>0$
\BQNY
\pk{u^2(T-\tau(u))> x|\tau(u)<T}
=
\frac{\pk{\tau(u)<T-xu^{-2}}}{\pk{\tau(u)<T}}.
\EQNY
With $T_u:=T- x u^{-2}$ and
$\widetilde{X}_u(t):=\frac{Y(t)}{\sigma(t)}\frac{g_u(T_u)}{g_u(t)}$ the above \rH{can be re-written as}
\begin{eqnarray*}
\pk{u^2(T-\tau(u))> x|\tau(u)<T}
=\frac{\pk{\sup_{t\in[0,T_u]}\widetilde{X}_u(t)>g_u(T_u)}}{\pk{\sup_{t\in[0,T]}X_u(t)>g_u(T)}}.
\end{eqnarray*}
As in the proof of Theorem \ref{ThmM} we have
\begin{eqnarray*}
\pk{\sup_{t\in[0,T_u]}\widetilde{X}_u(t)>g_u(T_u)}
\ge\Psi(g_u(T_u)).
\end{eqnarray*}
In order to derive the upper bound we use a time change such that
\BQNY
\pk{\sup_{t\in[0,T_u]}\widetilde{X}_u(t)>g_u(T_u)}=\pk{\sup_{t\in[0,1]}\widetilde{X}_u(T_u t)>g_u(T_u)}.
\EQNY
Similar argumentations as in \eqref{DDt} and \eqref{intruincovX} yield that, for some $\theta_0\in(0,1)$
$$
\sigma_{\widetilde{X}_u}(T_u t)\le \frac{1}{1+\frac{\sigma'(T_u)}{2\sigma(T_u)}T_u(1-t)}
$$
and
$$
r_{\widetilde{X}_u}(T_u s,T_u t)\ge e^{-2CT_u^2(t-s)^2}
$$
hold for all $s,t\in[\theta_0,1]$ and all $u$ sufficiently large. Consequently, in view of  iii) in \netheo{ThmPiter} and similar argumentations as in the proof of \netheo{ThmM} we conclude that
\begin{eqnarray*}
\pk{\sup_{t\in[0,T_u]}\widetilde{X}_u(t)>g_u(T_u)}
\le\Psi(g_u(T_u))(1+o(1))
\end{eqnarray*}
as $u\rw\IF.$
Hence
\begin{eqnarray}
\pk{u^2(T-\tau(u))> x|\tau(u)<T}
&=&
\frac{\Psi(g_u(T_u))}{\Psi(g_u(T))}(1+o(1))\nonumber\\
&=&
\exp\left(\frac{g_u^2(T)-g_u^2(T_u)}{2} \right)(1+o(1)), \quad u\to \IF.
\label{as1}
\end{eqnarray}
After some standard algebra, it follows that 
\begin{eqnarray}
g_u^2(T)-g_u^2(T_u)
=
\frac{(u+c\widetilde{\delta}(T))^2}{\sigma^2(T)}
-\frac{(u+c\widetilde{\delta}(T_u))^2}{\sigma^2(T_u)}
= -\frac{2\sigma'(T)}{\sigma^3(T)}x (1+o(1))
\end{eqnarray}
as $u\to\infty$. Consequently, by (\ref{as1})
\begin{eqnarray*}
\lim_{u\to\infty}
\pk{u^2(T-\tau(u))> x|\tau(u)<T}
=\exp\left(     - \frac{\sigma'(T)}{\sigma^3(T)}x    \right),
\end{eqnarray*}
which completes the proof.
\QED

\section{Appendix}

We give below an extension of Theorem D.3 in Piterbarg (1996) suitable for a family of Gaussian processes which is in \rH{particular useful} for the proof of our main results.
We first introduce two well-known constants appearing in the asymptotic theory of Gaussian processes. Let  $\{B_\alpha(t),t\ge0 \}$ be  a standard fractional Brownian motion  with Hurst index
$\alpha/2 \in (0,1]$ which is a centered Gaussian process with covariance  function
\BQNY
\mathbb{C}ov(B_\alpha(t),B_\alpha(s))=\frac{1}{2}(t^{\alpha}+s^{\alpha}-\mid t-s\mid^{\alpha}),\quad s,t\ge0.
\EQNY
The {\it Pickands constant} is defined by
\BQNY\label{pick}
\mathcal{H}_\alpha=\lim_{T\rightarrow\infty}\frac{1}{T} \E{ \exp\biggl(\sup_{t\in[0,T]}\Bigl(\sqrt{2}B_\alpha(t)-t^{\alpha}\Bigr)\biggr)} \cE{\in (0,\IF)}\ \alpha  \in (0,2]
 \EQNY
and the {\it Piterbarg constant} is given by
\BQNY
\mathcal{P}_\alpha^b=\lim_{T\rw\IF}\E{ \exp\biggl(\sup_{t\in[0,T]}\Bigl(\sqrt{2}B_\alpha(t)-(1+b)t^{\alpha}\Bigr)\biggr)} \in (0,\IF),\ \alpha  \in (0,2],\ b>0.
\EQNY
See for instance Piterbarg (1996) and D\c{e}bicki and Mandjes (2003) for properties of the above two constants.
Assume, in what follows, that $\theta$ and $ T$ are two positive constants satisfying $\theta<T$.
Let $\{\eta_u(t), t\ge0\}$ be a family of Gaussian processes satisfying the following three assumptions:

${\bf A1}:$ The variance function $\sigma_{\eta_u}^2(t)$ of $\eta_u$ attains its maximum over $[\theta,T]$ \rH{at} \Ke{the} unique point $t=T$ for any  $u$ large enough, and further there exist two positive constants $A, \beta$ and a function $A(u)$ satisfying $\lim_{u\rw\IF}A(u)=A$ such that $\sigma_{\eta_u}(t)$ has the following expansion around $T$ for all $u$ large enough
\BQNY
\sigma_{\eta_u} (t)=1-A(u)(T-t)^\beta(1+o(1)), \ \ t\uparrow T.
\EQNY

${\bf A2}:$ There exist two constants $ \alpha\in(0,2], B>0$ and a function $B(u)$ satisfying $\lim_{u\rw\IF}B(u)=B$ such that the correlation function $r_{\eta_u}(s,t)$ of $\eta_u$ has the following expansion around $T$ for all $u$ large enough
\BQNY
r_{\eta_u} (s,t)=1-B(u)\abs{t-s}^\alpha(1+o(1)), \ \ \min(s,t)\uparrow T.
\EQNY

${\bf A3}:$ For some positive constants $\mathbb{Q}$ and $\gamma$, and all $u$ large enough
$$
\E{\eta_u(s)-\eta_u(t)}^2\le \mathbb{Q}\abs{t-s}^\gamma
$$
for any $s,t\in[\theta,T]$.
\BT\label{ThmPiter}
Let $\{\eta_u(t), t\ge0\}$ be a family of Gaussian processes satisfying Assumptions {\bf A1-A3}.

i) If $\beta>\alpha$, then
$$
\pk{\sup_{t\in[\theta,T]}\eta_u(t)>u}= \frac{B^{\frac{1}{\alpha}}}{\sqrt{2\pi}A^{\frac{1}{\beta}}}
\mathcal{H}_{\alpha} \Gamma\left(\frac{1}{\beta}+1\right) u^{\frac{2}{\alpha}-\frac{2}{\beta}-1}
\exp\left(-\frac{u^2}{2}\right)\oo,\ \ \text{as}\ u\rw\IF.
$$
ii) For  $\beta=\alpha$ we have
$$
\pk{\sup_{t\in[\theta,T]}\eta_u(t)>u}=\frac{1}{\sqrt{2\pi}}\mathcal{P}_\alpha^{\frac{A}{B}}  u^{-1}\exp\left(-\frac{u^2}{2}\right)\oo,\ \ \text{as}\ u\rw\IF.
$$
iii) If $\beta<\alpha$, then
$$
\pk{\sup_{t\in[\theta,T]}\eta_u(t)>u}= \frac{1}{\sqrt{2\pi}} u^{ -1}\exp\left(-\frac{u^2}{2}\right)\oo,\ \ \text{as}\ u\rw\IF.
$$

\ET
\prooftheo{ThmPiter} 
Since
from assumptions {\bf A1-A2} we have that, for any $\vn>0$ and $u$ large enough
\BQNY
(A-\vn)(T-t)^\beta(1+o(1))\le 1-\sigma_{\eta_u} (t)\le (A+\vn)(T-t)^\beta(1+o(1)), \ \ t\uparrow T
\EQNY
and
\BQNY
(B-\vn)\abs{t-s}^\alpha(1+o(1))\le 1-r_{\eta_u} (s,t)\le (B+\vn)\abs{t-s}^\alpha(1+o(1)), \ \ \min(s,t)\uparrow T
\EQNY
Theorem D.3 in Piterbarg (1996) gives tight asymptotic upper and lower bounds, and thus the claims follow by \rH{letting} $\vn\rw0$. \QED

{\bf Acknowledgement}: K. D\c{e}bicki was partially supported by
NCN Grant No 2011/01/B/ST1/01521 (2011-2013).
All the authors kindly acknowledge partial support by the Swiss National Science Foundation Grant 200021-1401633/1 and
by the project RARE -318984, a Marie Curie IRSES Fellowship within the 7th
European Community Framework Programme.

\end{document}